\documentclass[12pt]{article}
\usepackage{amsfonts}
\usepackage{latexsym,amsmath,amssymb}
\usepackage{graphics,epstopdf}
\usepackage[pdftex]{graphicx}
\numberwithin{equation}{section}
\begin{document}
\bigskip
\bigskip
\begin{center}
{\large\bf{On the Fuzzy Stability of an Affine Functional Equation}}\\
\end{center}
\begin{center}
\vspace{.7cm}
{\bf{Md. Nasiruzzaman}} \\
\vspace{.2cm}
Department of Mathematics, Aligarh Muslim University, Aligarh 202002, India\\
Email: nasir3489@gmail.com\\
\end{center}
{\large\bf{Abstract:}} In this paper, we obtain the general solution of the following functional equation
$$f(3x+y+z)+f(x+3y+z)+f(x+y+3z)
+f(x)+f(y)+f(z)=6f(x+y+z).$$
We establish the Hyers-Ulam-Rassias stability of the above functional equation in the fuzzy normed spaces.
Further we show the above functional equation is stable in the sense of Hyers and Ulam in fuzzy normed spaces.
 \vspace{.2cm}\\
{\large\bf{1. Introduction}}
\vspace{.2cm}\\
$~~~~~$In modelling applied problems only partial informations may be known (or) there may be a degree of uncertainty in the
 parameters used in the model or some measurements may be imprecise. Due to such features, we are tempted to consider
 the study of functional equations in the fuzzy setting.\\
For the last 40 years, fuzzy theory has become very active area of research and a lot of development has been made in
 the theory of fuzzy sets \cite{Zadeh} to find the fuzzy analogues of the classical set theory. This branch finds a wide range of
 applications in the field of science and engineering.\\
A.K. Katsaras \cite{Katsaras} introduced an idea of fuzzy norm on a linear space in 1984, in the same year Cpmgxin Wu and Jinxuan
Fang \cite{Fang} introduced a notion of fuzzy normed space to give a generalization of the Kolmogoroff normalized theorem
 for fuzzy topological linear spaces. In 1991, R. Biswas \cite{Biswas} defined and studied fuzzy inner product spaces in
 linear space. In 1992, C. Felbin \cite{Felbin} introduced an alternative definition of a fuzzy norm on a linear topological
  structure of a fuzzy normed linear spaces. In 2003, T. Bag and S.K. Samanta \cite{Bag1} modified the definition of S.C. Cheng
   and J.N. Mordeson \cite{Mordeson} by removing a regular condition. \\
In 1940, Ulam \cite{ulam} raised a question concerning the stability of group homomorphism as follows:\\
Let $G_1$ be a group and $G_2$ a metric group with the metric $d(.,.).$ Given $\varepsilon>0$, does there exists a $\delta >0$ such that if a function $f:G_1\to G_2$ satisfies the inequality $$d(f(xy),f(x)f(y))<\delta ~\mbox{for all}~ x,y\in G_1,$$ then there exists a homomorphism $h:G_1\to G_2$ with $$d(f(x),H(x))<\varepsilon ~\mbox{for all}~ x\in G_1?$$ The concept of stability for a functional equation arises when we replace the functional equation by an inequality which acts as a perturbation of the equation. 
In 1941, the case of approximately additive mappings was solved by Hyers \cite{hyers} under the assumption that $G_2$ is a Banach space. In 1978, a generalized version of the theorem of Hyers for approximately linear
mapping was given by Th.M. Rassias \cite{rassias3}. He proved that for a mapping $f:E_1\to E_2$ such that $f(tx)$ is continuous in $t\in \mathbb{R}$ and for each fixed $x\in E_1$ assume that
there exist a constant $\varepsilon>0$ and $p\in[0,1)$ with
$$\parallel f(x + y)-f(x)-f(y)\parallel\leqslant \varepsilon(\parallel x\parallel^p+\parallel y\parallel^p)\eqno(1.1)$$
$x,y\in E_1$, then there exist a unique $R$-Linear mapping $T: E_1\to E_2$ such that
$$\parallel f(x)-T(x)\parallel\leqslant \frac{2\varepsilon}{2-2^p}\parallel x\parallel^p~~~(x\in E_1)\eqno(1.2)$$
\\
The result of Rassias has influenced the development of what is now called the Hyers-Ulam-Rassias stability theory
for functional equations. In 1994, a generalization of Rassias theorem was obtained by Gavruta \cite{gavruta1} by replacing the
bound $\varepsilon(\|x\|^p+\|y\|^p)$ by a general control function $\varphi(x, y)$.
During the last decades, the stability problems of several functional equations
have been extensively investigated by a number of authors (c.f. \cite{czerwik}, \cite{hyers1}, \cite{kannappan}, \cite{ansari} and \cite{sam1}--\cite{mohimcm} etc.).
In 1982-1989, J.M.Rassias \cite{jmrassias1,jmrassias2} replaced the sum appeared in right hand side
of the equation (1.1) by the product of powers of norms. In fact, he proved the
following theorem.\\
\\
{\large\bf{Theorem 1.1}} Let $f: E_1\to E_2$ be a mapping from a normed vector space $E_1$ into
Banach space $E_2$ subject to the inequality
$$\parallel f(x + y)-f(x)-f(y)\parallel\leqslant \varepsilon(\parallel x\parallel^p\parallel y\parallel^p)\eqno(1.3)$$
for all $x,y\in E_1$, where $\varepsilon$ and $p$ are constants with $\varepsilon > 0$ and $0 \leqslant p < \frac12$. Then the limit
$$L(x)=\lim\limits_{n\to\infty}\frac{f(2^nx)}{2^n}\eqno(1.4)$$
exists for all $x\in E_1$, and $L: E_1\to E_2$ is the unique additive mapping which satisfies
$$\parallel f(x)-L(x)\parallel\leqslant \frac{\varepsilon}{2-2^{2p}}\parallel x\parallel^{2p}\eqno(1.5)$$
for all $x\in E_1$. If $p > \frac12$ the inequality (1.3) holds for $x,y \in E_1 $ and the limit
$$A(x)=\lim\limits_{n\to\infty}2^nf\bigl{(}\frac{x}{2^n}\bigl{)}\eqno(1.6)$$
exists for all $x\in E_1$ and $A : E_1\to E_2$ is the unique additive mapping which satisfies
$$\parallel f(x)-A(x)\parallel\leqslant \frac{\varepsilon}{2^{2p}-2}\parallel x\parallel^{2p}~~~~(\forall x\in E_1)\eqno(1.7)$$
\\
 Recently, Cadariu et al \cite{cadariu} studied the generalized Hyers-Ulam stability by using the direct method as well as the fixed point method for the affine type functional equation
$$f(2x+y)+f(x+2y)+f(x)+f(y)=4f(x+y), ~~\mbox{for all}~ x,y\in G.\eqno(1.8)$$
In the present paper, we obtain the general solution of the following functional equation
$$f(3x+y+z)+f(x+3y+z)+f(x+y+3z)
+f(x)+f(y)+f(z)=6f(x+y+z).\eqno(1.9)$$
where $f:X\to Y$, $X$ and $Y$ are normed spaces. Then, we establish the fuzzy Hyers-Ulam-Rassias stability of the above functional equation.\\
\\
{\large\bf{2. Preliminary Notes}}\\
\\
Before we proceed to the main results, we will introduce a definition and
some examples to illustrate the idea of fuzzy norm.
\vspace{.2cm}\\
\\
{\large\bf{Definition 2.1}} Let X be a real linear space. A mapping $N:X\times\mathbb{R}\to[0,1]$
 (the so-called fuzzy ~subset) is said to be a $fuzzy ~norm$ on $X$ if for all $x,y\in X$ and all $s,~t\in \mathbb{R}$,\\
$(N_1)~N(x,t)=0~\mbox{for}~t\leqslant0;$\\
$(N_2)~x=0~\mbox{if and only if}~N(x,t)=1~\mbox{for all}~t>0$;\\
$(N_3)~N(cx,t)=N(x,t/\mid c\mid)~\mbox{if}~c\neq0$;\\
$(N_4)~N(x+y,t+s)\geqslant\min\{N(x,t),N(y,s)\}$;\\
$(N_5)~N(x,.)~\mbox{is a non-decreasing function on $\mathbb{R}$ and}~\lim\limits_{t\to\infty}N(x,t)=1$;\\
$(N_6)~\mbox{for}~x\neq0,~N(x,.)~\mbox{is continuous on $\mathbb{R}$}$.\\
The pair $(X,N)$ is called a $fuzzy ~normed ~linear ~space$. One may regard $N(x,t)$ as the truth value of
 the statement that the norm of $x$ is less than or equal to the real number $t$.
 \vspace{.2cm}\\
 \\
{\large\bf{Example 2.2}} Let $(X,\|.\|)$ be a normed linear space. One can be easily verify that for each
 $p>0$,
 $$ N_p(x,t)={\begin{cases} \frac{t}{t+p\|x\|} ~~~t>0, x\in X\cr$$
$$ 0 ~~~~~~~~~t\leqslant0, x\in X\cr\end{cases}}$$
 is a fuzzy norm on $X$.
 \vspace{.2cm}\\
 \\
{\large\bf{Example 2.3}} Let $(X,\|.\|)$ be a normed linear space. The mapping
 $N: X\times\mathbb{R}\to[0,1]$ by\\
$$N(x,t)= {\begin{cases} \frac{t^2-\|x\|^2}{t^2+\|x\|^2} ~~~t>\|x\| \cr$$
$$ 0 ~~~~~~~~~~t\leqslant\|x\| \cr\end{cases}}$$\\ is a fuzzy norm on $X$.
\vspace{.2cm}\\
\\
{\large\bf{Definition 2.4}} Let $(X,N)$ be a fuzzy normed linear space. A sequence $\{x_n\}$ in $X$ is said to
be convergent if there exists an $x\in X$ such that $\lim\limits_{n\to\infty}N(x_n-x,t)=1$
 for all $t>0$. In this case, $x$ is called the limit of the sequence $\{x_n\}$ and we denote it by
$$N-\lim\limits_{n\to\infty}N(x_n-x,t)=x.$$
\\
{\large\bf{Definition 2.5}} Let $(X,N)$ be a fuzzy normed linear space. A sequence $\{x_n\}$ in $X$ is said to be
 Cauchy if for each $\varepsilon>0$ and each $\delta>0$ there exists an $n_0\in \mathbb{N}$ such that
  $$N(x_m-x_n,\delta)>1-\varepsilon~~(m,n\geqslant n_0).$$
  It is well known that every convergent sequence in a fuzzy normed linear space is Cauchy. If each Cauchy
  sequence is convergent, then the fuzzy norm is said to be complete and the fuzzy normed vector space is called  a fuzzy Banach space.\\
  \\
The remaining part of the paper is organized as follows: We discuss the general solution of functional equation (1.9) in Section 3. Section 4 is devoted to investigate the non-uniform version of stability of functional equation (1.9) in fuzzy normed spaces and in section (5), we show under suitable conditions that in fuzzy normed spaces functional equation (1.9) is stable uniformly. \\
Now we proceed to find the general solution of the functional equation (1.9)
\vspace{.2cm}\\
{\large\bf{3. Solution of the Functional Equation (1.9)}}
\vspace{.2cm}\\
{\large\bf{Theorem 3.1}} A mapping $f:X\to Y$, $X$ and $Y$ are normed spaces,
is a solution of the functional equation (1.9) if and only if it is an affine mapping
(i.e., it is the sum between a constant and an additive function).
\vspace{.2cm}\\
{\large\bf{Proof.}} We can easily seen that any affine function $f$ is a solution of the equation (1.9).\\
Conversely, we have two cases:\\
$Case ~1:$ $f(0)=0$.\\
If we take $y=z=-x$ in (1.9), we obtain
$$2f(x)+2f(-3x)+2f(-x)=6f(-x),~~\mbox{for all}~ x\in X.\eqno(3.1)$$
Again replacing putting $y=z=0$ in (1.9), we obtain
$$f(3x)=3f(x),~~\mbox{for all} ~x\in X.\eqno(3.2)$$
By (3.1) and (3.2), we have
$f(-x)=-f(x),~\mbox{for all}~x\in X $. It results that $f$ is an odd mapping.
Replace $z$ by $-y$ in (1.9), we get
$$f(x+2y)+f(x-2y)=2f(x)\eqno(3.3)$$
If we replace $x$ and $y$ by $\frac{u+v}2$ and $\frac{u-v}{4}$, respectively, in (3.3) and using (3.2),
we have $$f(u+v)=f(u)+f(v),~~\mbox{for all}~ u,v\in X.$$
So, $f$ is an additive mapping.\\
$Case ~2:$ General case.
Let us consider the function $g(x):=f(x)-f(0)$. It is clear that $g(0)=0$ and $f(x)=g(x)+f(0)$.\\
Replacing $f$ by $g$ in (1.9), it results
$$g(3x+y+z)+g(x+3y+z)+g(x+y+3z)
+g(x)+g(y)+g(z)=6g(x+y+z).$$
$\mbox{for all}~x,y,z\in X$.
Taking in account that $g(0)=0$, from $Case~1$, we obtain that $g$ is an additive mapping, hence \\
$f(x)=g(x)+f(0)$
is an affine function.\\
This completes the proof.
\vspace{.2cm}\\
\vspace{.2cm}\\
For a given mapping $f:X\to Y$, let us denote
\begin{align*}
Df(x,y,z)&=f(3x+y+z)+f(x+3y+z)+f(x+y+3z)\\
&+f(x)+f(y)+f(z)-6f(x+y+z)
\end{align*}
{\large\bf{4. Fuzzy Hyers-Ulam-Rassias Stability: non-uniform version}}\\
\\
{\large\bf{Theorem 4.1}} Let $X$ be a linear space and $(Z,N')$ a fuzzy normed space. Let\\
 $\varphi:X^3\to Z$
be a mapping such that for some $\alpha\neq0$ with $0<\alpha< 3$
$$N'(\varphi(3x,0,0),t)\geqslant N'(\alpha\varphi(x,0,0),t)\eqno(4.1)$$
for all $x\in X$, $t>0$ and $$\lim\limits_{n\to\infty}N'(\varphi(3^nx,3^ny,3^nz), 3^nt)=1,$$
for all $x,y,z\in X$ and all $t>0$. Suppose that $(Y,N)$ be a fuzzy Banach space and an odd mapping $f: X\to Y$ satisfies the inequality
$$N(Df(x,y,z),t)\geqslant N'(\varphi(x,y,z), t)\eqno(4.2)$$ for all $x,y,z\in X$
and all $t>0$.
Then the limit $$A(x)=N-\lim\limits_{n\to\infty}\frac{f(3^nx)}{3^n}$$ exists for all $x\in X$ and the mapping
$A: X\to Y$ is the unique affine mapping satisfying $$N(f(x)-A(x)-f(0), t)\geqslant N'(\varphi(x,0,0),
(3-\alpha)t)\eqno(4.3)$$
for all $x\in X$ and all $t>0$.
\vspace{.2cm}\\
{\large\bf{Proof.}} Letting $y=z=0$ in (4.2), we get
$$N(f(3x)-3f(x)+2f(0), t)\geqslant N'(\varphi(x,0,0), t)\eqno(4.4)$$ for all $x\in X$ and all $t>0$.\\
If we define the mapping $g: X\to Y$ such that $g(x):=f(x)-f(0)$ for all $x\in X$. Indeed $g(0)=0$.
Then (4.4) implies
$$N(g(3x)-3g(x), t)\geqslant N'(\varphi(x,0,0), t)$$
Replacing $x$ by $3^nx$ in the last inequality, we obtain
$$N(g(3^{n+1}x)-3g(3^nx), t)\geqslant N'(\varphi(3^nx,0,0), t)$$
$$N\biggl{(}\frac{g(3^{n+1}x)}{3^{n+1}}-\frac{g(3^nx)}{3^n}, \frac t{3^{n+1}}\biggl{)}\geqslant N'(\varphi(x,0,0), \frac{t}{\alpha^n})$$
$$N\biggl{(}\frac{g(3^{n+1}x)}{3^{n+1}}-\frac{g(3^nx)}{3^n}, \frac {\alpha^nt}{3^{n+1}}\biggl{)}\geqslant N'(\varphi(x,0,0), t)\eqno(4.5)$$
for all $x\in X$ and all $t>0$.
It follows from $\frac{g(3^nx)}{3^n}-g(x)=\sum\limits_{j=0}^{n-1}\frac{g(3^{j+1}x)}{3^{j+1}}-\frac{g(3^jx)}{3^j}$
and (4.5) that
\begin{align*}
    N\biggl{(}\frac{g(3^nx)}{3^n}-g(x),\sum\limits_{j=0}^{n-1}\frac{\alpha^jt}{3^{j+1}}\biggl{)}&=N\biggl{(}\sum\limits_{j=0}^{n-1}\frac{g(3^{j+1}x)}{3^{j+1}}-\frac{g(3^jx)}{3^j}, \sum\limits_{j=0}^{n-1}\frac {\alpha^jt}{3^{j+1}})\biggl{)}\\
    &\geqslant\min\bigcup\limits_{j=0}^{n-1}\biggl{\{}N\biggl{(}\frac{g(3^{j+1}x)}{3^{j+1}}-\frac{g(3^jx)}{3^j}, \frac {\alpha^jt}{3^{j+1}}\biggl{)}\biggl{\}}\\
    &\geqslant N'(\varphi(x,0,0),t).............................................(4.6)
\end{align*}
for all $x\in X$ and all $t>0$. Replacing $x$ by $3^mx$ in (4.6), we get
$$N\biggl{(}\frac{g(3^{n+m}x)}{3^{n+m}}-\frac{g(3^mx)}{3^m}, \sum\limits_{j=0}^{n-1}\frac {\alpha^jt}{3^{j+m+1}})\biggl{)}\geqslant N'\biggl{(}\varphi(x,0,0),\frac t{\alpha^m}\biggl{)}$$
and so $$N\biggl{(}\frac{g(3^{n+m}x)}{3^{n+m}}-\frac{g(3^mx)}{3^m}, \sum\limits_{j=m}^{n+m-1}\frac {\alpha^jt}{3^{j+1}})\biggl{)}\geqslant N'(\varphi(x,0,0),t)$$
$$N\biggl{(}\frac{g(3^{n+m}x)}{3^{n+m}}-\frac{g(3^mx)}{3^m}, t\biggl{)}\geqslant N'\biggl{(}\varphi(x,0,0),\frac t{\sum\limits_{j=m}^{n+m-1}\frac {\alpha^j}{3^{j+1}}}\biggl{)}\eqno(4.7)$$
for all $x\in X$, $t>0$ and $m,n\geqslant0$. Since $0<\alpha< 3 $ and $\sum\limits_{j=0}^{\infty}(\frac\alpha 3)^j<\infty$,
 the Cauchy criterion for convergence and $(N_5)$ imply that $\{\frac{g(3^nx)}{3^n}\}$ is a Cauchy sequence in $(Y,N)$.
 Since $(Y,N)$ is a fuzzy Banach space, this sequence converges to some point $A(x)\in Y$. Hence, we can define a mapping
 $A: X\to Y$ by $A(x)=N-\lim\limits_{n\to\infty}\frac{g(3^nx)}{3^n}=N-\lim\limits_{n\to\infty}\frac{f(3^nx)}{3^n}$ for all $x\in X$, namely.
 Since $f$ is odd, $A$ is odd. Letting $m=0$ in (4.7), we get
 $$N\biggl{(}\frac{g(3^{n}x)}{3^{n}}-g(x), t\biggl{)}\geqslant N'\biggl{(}\varphi(x,0,0),\frac t{\sum\limits_{j=0}^{n-1}\frac {\alpha^j}{3^{j+1}}}\biggl{)}$$
 Taking the limit as $n\to\infty$ and using $(N_6)$, we get
 \begin{align*}
    N(A(x)-g(x),t)&\geqslant N'\biggl{(}\varphi(x,0,0),\frac t{\sum\limits_{j=0}^{\infty}\frac {\alpha^j}{3^{j+1}}}\biggl{)}\\
    &=N'(\varphi(x,0,0),(3-\alpha)t)
 \end{align*}
$$N(f(x)-A(x)-f(0),t)\geqslant N'(\varphi(x,0,0),(3-\alpha)t)$$ for all $x\in X$ and all $t>0$.\\
Now we claim that $A$ is affine. Replacing $x,y,z$ by $3^nx,3^ny,3^nz$, respectively, in (4.2), we get
$$N\biggl{(}\frac1{3^n}Df(3^nx,3^ny,3^nz),t\biggl{)}\geqslant N'(\varphi(3^nx,3^ny,3^nz),3^nt)$$
for all $x,y,z\in X$ and all $t>0$. Since
$$\lim\limits_{n\to\infty}N'(\varphi(3^nx,3^ny,3^nz), 3^nt)=1,$$ $A$ satisfies the functional equation (1.9). Hence $A$ is affine. To prove the uniqueness of $A$, let $A': X\to Y$ be another affine mapping satisfying (4.3). Fix $x\in X$. Clearly $A(3^nx)=3^nA(x)$ and $A'(3^nx)=3^nA'(x)$ for all $x\in X$ and all $n\in \mathbb{N}$. It follows from (4.3) that
\begin{align*}
    N(A(x)-A'(x),t)&=N\biggl{(}\frac{A(3^nx)}{3^n}-\frac{A'(3^nx)}{3^n},t\biggl{)}\\
    &\geqslant\min\biggl{\{}N\biggl{(}\frac{A(3^nx)}{3^n}-\frac{g(3^nx)}{3^n},\frac t2\biggl{)},~N\biggl{(}\frac{g(3^nx)}{3^n}-\frac{A'(3^nx)}{3^n},\frac t2\biggl{)}\biggl{\}}\\
    &\geqslant N'\biggl{(}\varphi(3^nx,0,0),\frac{3^n(3-\alpha)t}{2}\biggl{)}\\
    &\geqslant N'\biggl{(}\varphi(x,0,0),\frac{3^n(3-\alpha)t}{2\alpha^n}\biggl{)}
\end{align*}
for all $x\in X$ and all $t>0$. Since $\lim\limits_{n\to\infty}\frac{3^n(3-\alpha)}{2\alpha^n}=\infty$,
we obtain
$$\lim\limits_{n\to\infty} N'\biggl{(}\varphi(x,0,0),\frac{3^n(3-\alpha)t}{2\alpha^n}\biggl{)}=1.$$
Thus $N(A(x)-A'(x),~t)=1$ for all $x\in X$ and all $t>0$, and so $A(x)=A'(x)$.\\
This completes the proof.
\vspace{.2cm}\\
{\large\bf{5. Fuzzy Hyers-Ulam-Rassias Stability: uniform version}}\\
\\
{\large\bf{Theorem 5.1}} Let $X$ be a linear space and $(Y,N)$ be a fuzzy Banach space. Let\\
$\varphi:~X^3\to[0,\infty)$ be a function such that
$$\tilde{\varphi}(x,y,z)=\sum\limits_{n=0}^{\infty}\frac1{3^n}\varphi(3^nx,3^ny,3^nz)<\infty\eqno(5.1)$$
for all $x,y,z\in X$. Let $f:X\to Y$ be a uniformly approximately affine mapping with respect to
$\varphi$ in the sense that $$\lim\limits_{t\to\infty}N(Df(x,y,z),t\varphi(x,y,z))=1\eqno(5.2)$$
uniformly on $X^3$. Then
$$A(x):= N-\lim\limits_{n\to\infty}\frac{f(3^nx)}{3^n}$$ for all $x\in X$ exists and defines an affine mapping $A:~X\to Y$ such that if for some $~\alpha>0,~\delta>0$
$$N(Df(x,y,z),\delta\varphi(x,y,z))>\alpha\eqno(5.3)$$ for all $x,y,z\in X$, then
$$N(f(x)-A(x)-f(0),\frac{\delta}3\tilde{\varphi}(0,0,,x))>\alpha$$ for all $x\in X$.
\vspace{.2cm}\\
{\large\bf{Proof.}} Let $\varepsilon>0$, by (5.2), we can find $t_0>0$ such that
$$N(Df(x,y,z),t\varphi(x,y,z))\geqslant1-\varepsilon\eqno(5.4)$$ for all $x,y,z\in X$ and all $t\geqslant t_0$. Define $g:~X\to Y$ such that $g(x):=f(x)-f(0)$. It is clear that $g(0)=0$ and $f(x)=g(x)+f(0)$. Now $(5.4)$ implies that
$$N(Dg(x,y,z),t\varphi(x,y,z))\geqslant1-\varepsilon\eqno(5.5)$$ for all $x,y,z\in X$
and all $t\geqslant t_0$. By induction on $n$, we will show that
$$N\biggl{(}g(3^nx)-3^ng(x),t\sum\limits_{m=0}^{n-1}3^{n-m-1}\varphi(0,0,3^mx)\biggl{)}\geqslant1-\varepsilon\eqno(5.6)$$ for all $x\in X$, all $t\geqslant t_0$ and $n\in \mathbb{N}$. Putting $x=y=0$ and $z=x$ in (5.5), we get (5.6) for $n=1$. Let (5.6) holds for some positive integers $n$. Then
   $$ N(g(3^{n+1}x)-3^{n+1}g(x),t\sum\limits_{m=0}^{n}3^{n-m}\varphi(0,0,3^mx))$$
    $$\geqslant\min\{N(g(3^{n+1}x)-3g(3^nx),t\varphi(0,0,3^nx)),$$
   $$ ~N(3g(3^nx)-3^{n+1}g(x),t\sum\limits_{m=0}^{n}3^{(n-m)}\varphi(0,0,3^mx))\}$$
     $$\geqslant\min\{1-\varepsilon,1-\varepsilon\}=1-\varepsilon.$$
This completes the induction argument. Let $t=t_0$ and put $n=p$. Then by replacing $x$ with $3^nx$ in (5.6), we obtain
$$N(g(3^{n+p}x)-3^pg(3^nx),t_0\sum\limits_{m=0}^{p-1}3^{p-m-1}\varphi(0,0,3^{n+m}x))\geqslant1-\varepsilon$$
$$N\biggl{(}\frac{g(3^{n+p}x)}{3^{n+p}}-\frac{g(3^nx)}{3^n},t_0\sum\limits_{m=0}^{p-1}3^{-(n+m+1)}\varphi(0,0,3^{n+m}x)\biggl{)}\geqslant1-\varepsilon\eqno(5.7)$$
for all integers $n\geqslant0,~p>0$. The convergence of (5.1) and the equation
$$\sum\limits_{m=0}^{p-1}3^{-(n+m+1)}\varphi(0,0,3^{n+m}x))=\frac1k\sum\limits_{m=n}^{n+p-1}3^{-m}\varphi(0,0,3^mx)$$
guarantees that for given $\delta>0$, there exists $n_0\in\mathbb{N}$ such that
$$\frac{t_0}3\sum\limits_{m=n}^{n+p-1}3^{-m}\varphi(0,0,3^mx)<\delta$$ for all $n\geqslant n_0$ and $p>0$. It follows from (5.7) that
$$N\biggl{(}\frac{g(3^{n+p}x)}{3^{n+p}}-\frac{g(3^nx)}{3^n},~\delta\biggl{)}\geqslant N\biggl{(}\frac{g(3^{n+p}x)}{3^{n+p}}-\frac{g(3^nx)}{3^n},~t_0\sum\limits_{m=0}^{p-1}3^{-(n+m+1)}\varphi(0,0,3^{n+m}x)\biggl{)}\geqslant1-\varepsilon\eqno(5.8)$$
for each $n\geqslant n_0$ and all $p>0$. Hence $\{\frac{g(3^nx)}{3^n}\}$ is a Cauchy sequence in $Y$. Since $Y$ is a fuzzy Banach space, this sequence converges to some $A(x)\in Y$. Hence we can define a mapping $A:X\to Y$ by
$A(x):=N-\lim\limits_{n\to\infty}\frac{g(3^nx)}{3^n}=N-\lim\limits_{n\to\infty}\frac{f(3^nx)}{3^n}$, for all $x\in X$ namely. For each $t>0$ and $x\in X$ $$\lim\limits_{n\to\infty}N\biggl{(}A(x)-\frac{f(3^nx)}{3^n},t\biggl{)}=1.$$
Now, let $x,y,z\in X$. Fix $t>0$ and $0<\varepsilon<1$. Since $\lim\limits_{n\to\infty}\frac1{3^n}\varphi(3^nx,3^ny,3^nz)=0$, there is some $n_1>n_0$ such that
$$N(DA(x,y,z),t)\geqslant\min\biggl{\{}N\biggl{(}A(3x+y+z)-\frac{f(3^n(3x+y+z))}{3^n},\frac{t}{8}\biggl{)},$$
$$N\biggl{(}A(x+3y+z)-\frac{f(3^n(x+3y+z))}{3^n},\frac{t}{8}\biggl{)},$$
$$N\biggl{(}A(x+y+3z)-\frac{f(3^n(x+y+3z))}{3^n},\frac{t}{8}\biggl{)},$$
$$N\biggl{(}A(x)-\frac{f(3^nx)}{3^n},\frac{t}{8}\biggl{)},~N\biggl{(}A(y)-\frac{f(3^ny)}{3^n},\frac{t}{8}\biggl{)}, N\biggl{(}A(z)-\frac{f(3^nz)}{3^n},\frac{t}{8}\biggl{)},$$
$$~N\biggl{(}A(x+y+z)-\frac{f(3^n(x+y+z))}{3^n},\frac{t}{(8)}\biggl{)},$$
$$N\biggl{(}Df(3^nx,3^ny,3^nz),\frac{3^nt}{8}\biggl{)}\biggl{\}}$$
     The first $7$ terms on the right hand side of the above inequality tend to 1 as $n\to\infty$ and the last term is greater than $N(Df(3^nx,3^ny,3^nz),t_0\varphi(3^nx,3^ny,3^nz)),$ i.e., by (5.4), greater than or equal to $1-\varepsilon$. Thus $N(DA(x,y,z),t)\geqslant1-\varepsilon$ for all $t\geqslant0$ and $0<\varepsilon<1$.
It follows that $N(DA(x,y,z),t)=1$ for all $t>0$ and by $(N_2)$, we have
$DA(x,y,z)=1$, i.e.,
  $$A(3x+y+z)+A(x+3y+z)+A(x+y+3z)
+A(x)+A(y)+A(z)=6A(x+y+z)$$
To end the proof, let for some positive $\alpha$ and $\delta$, (5.3) holds. Let
$$\varphi_n(x,y,z):=\sum\limits_{m=0}^{n-1}3^{-(m+1)}\varphi(3^mx,3^my,3^mz)$$ for all $x,y,z\in X$. Let $x\in X$. By a similar discussion as in the begining of the proof, we can obtain from (5.3) $$N(g(3^nx)-3^ng(x),\delta\sum\limits_{m=0}^{n-1}3^{(n-m-1)}\varphi_n(0,0,3^mx))\geqslant\delta\eqno(5.9)$$
for all $n\in \mathbb{N}$. Let $s>0$. We have
$$N(g(x)-A(x),\delta\varphi_n(0,0,x)+s)\geqslant\min\biggl{\{}N\biggl{(}g(x)-\frac{g(3^nx)}{3^n},\\
~\delta\varphi_n(0,0,x)\biggl{)},~N\biggl{(}\frac{g(3^nx)}{3^n}-A(x),s\biggl{)}\biggl{\}}\eqno(5.10)$$
Combining (5.8), (5.9) and the fact that
$$\lim\limits_{n\to\infty}N\biggl{(}\frac{g(3^nx)}{3^n}-A(x),s\biggl{)}\\
=\lim\limits_{n\to\infty}N\biggl{(}\frac{f(3^nx)}{3^n}-A(x),s\biggl{)}=1,$$
we obtain that $$N(g(x)-A(x),\delta\varphi_n(0,0,x)+s)\geqslant\alpha$$ for large enough $n$. By the (upper semi) continuity of real function $N(g(x)-A(x),.)$, we obtain that
$$N\biggl{(}g(x)-A(x),\frac{\delta}3\tilde{\varphi}(0,0,x)+s\biggl{)}\geqslant\alpha.$$
Taking the limit as $s\to0$, we conclude that
$$N\biggl{(}g(x)-A(x),\frac{\delta}3\tilde{\varphi}(0,0,x)\biggl{)}\geqslant\alpha$$
$$N\biggl{(}f(x)-A(x)-f(0),\frac{\delta}3\tilde{\varphi}(0,0,x)\biggl{)}\geqslant\alpha.$$
This completes the proof.
\vspace{.2cm}\\
\\
{\large\bf{Theorem 5.2}} Let $X$ be a linear space and $(Y,N)$ be a fuzzy Banach space. Let\\
$\varphi:~X^3\to[0,\infty)$ be a function satisfying (5.1). Let
$f:~X\to Y$ be a uniformly approximately affine mapping with respect to $\varphi$. Then there is a unique affine mapping $A:~X\to Y$ such that
$$\lim\limits_{t\to\infty}N(f(x)-A(x)-f(0),t\tilde{\varphi}(0,0,x))=1\eqno(5.11)$$ uniformly on $X$.
\vspace{.2cm}\\
{\large\bf{Proof.}} The existence of uniform limit (5.11) immediately follows from Theorem 4.5. It remains to prove the uniqueness assertion. Let $A'$ be another affine mapping satisfying (5.11). Fix $c>0$. Given $\varepsilon>0$, by (5.11) for $A$ and $A'$, we can find some $t_0>0$ such that
$$N(g(x)-A(x),\frac t2\tilde{\varphi}(0,0,x))\geqslant1-\varepsilon,$$
$$N(g(x)-A'(x),\frac t2\tilde{\varphi}(0,0,x))\geqslant1-\varepsilon$$ for all $x\in X$ and all $t\geqslant t_0$.
Fix some $x\in X$ and find some integer $n_0$ such that\\
$t_0\sum\limits_{m=n}^{\infty}3^{-m}\varphi(0,0,3^mx)<\frac c2$, for all $n\geqslant n_0.$
Since
\begin{align*}
    \sum\limits_{m=n}^{\infty}3^{-m}\varphi(0,0,3^mx)&=\frac1{3^n}\sum\limits_{m=n}^{\infty}3^{-(m-n)}\varphi(0,0,3^{m-n}(3^nx))\\
    &=\frac1{3^n}\sum\limits_{j=0}^{\infty}\frac1{3^j}\varphi(0,0,3^j(3^nx))\\
    &=\frac1{3^n}\tilde{\varphi}(0,0,3^nx)
\end{align*}
We have
    $$N(A'(x)-A(x),c)\geqslant\min\biggl{\{}N\biggl{(}\frac{g(3^nx)}{3^n}-A(x),\frac c2\biggl{)},
    ~N\biggl{(}A'(x)-\frac{g(3^nx)}{3^n},\frac c2\biggl{)}\biggl{\}}$$
    $$=\min\biggl{\{}N\biggl{(}\frac{g(3^nx)}{3^n}-\frac{A(3^nx)}{3^n},\frac c2\biggl{)},
    ~N\biggl{(}\frac{A'(3^nx)}{3^n}-\frac{g(3^nx)}{3^n},\frac c2\biggl{)}\biggl{\}}$$
    $$=\min\biggl{\{}N\biggl{(}g(3^nx)-A(3^nx),\frac {3^nc}2\biggl{)},
    ~N\biggl{(}A'(3^nx)-g(3^nx),\frac{3^nc}2\biggl{)}\biggl{\}}$$
    $$\geqslant\min\biggl{\{}N\biggl{(}g(3^nx)-A(3^nx),3^nt_0\sum\limits_{m=n}^{\infty}3^{-m}\varphi(0,0,3^mx)\biggl{)},$$
    $$N\biggl{(}A'(3^nx)-g(3^nx),3^nt_0\sum\limits_{m=n}^{\infty}k^{-m}\varphi(0,0,3^mx)\biggl{)}\biggl{\}}$$
    $$=\min\biggl{\{}N\biggl{(}g(3^nx)-A(3^nx),t_0\tilde{\varphi}(0,0,3^nx)\biggl{)},$$
    $$N\biggl{(}A'(3^nx)-g(3^nx),t_0\tilde{\varphi}(0,0,3^mx)\biggl{)}\biggl{\}}$$
    $$\geqslant1-\varepsilon.$$
It follows that $N(A'(x)-A(x),~c)=1$, for all $c>0$. Thus $A(x)=A'(x)$ for all $x\in X$.\\
This completes the proof.\\
\\
Considering the control function $\varphi(x,y,z)=\varepsilon(\|x\|^p+\|y\|^p+\|z\|^p)$ for some $\varepsilon>0$, we obtain the following:\\
\\
{\large\bf{Corollary 5.3}} Let $X$ be a normed linear space, let $(Y,N)$ be a fuzzy Banach space, let $\varepsilon\geqslant0$, and let $0\leqslant p<1$. Suppose that $f:X\to Y$ is a function such that \\
$$\lim\limits_{n\to\infty}N(Df(x,y,z),~t\varepsilon(\|x\|^p+\|y\|^p+\|z\|^p))=1$$
uniformly on $X^3$. Then there is a unique affine mapping $A:X\to Y$ such that
$$\lim\limits_{t\to\infty}N\biggl{(}f(x)-A(x)-f(0),\frac{\varepsilon t3^{1-p}\|x\|^p}{3^{1-p}-1}\biggl{)}=1$$
uniformly on X.\\

\end{document}